\documentclass{article}
\usepackage{amsmath}
\usepackage{amssymb}

\newtheorem{theorem}{Theorem}

\newtheorem{definition}[theorem]{Definition}

\newtheorem{proposition}[theorem]{Proposition}

\input{tcilatex}
\begin{document}

\title{Klein geometries of high order}
\author{Erc\"{u}ment H. Orta\c{c}gil}
\maketitle

\begin{abstract}
Using irreducable representations of semi simple Lie algebras, we construct
Klein pairs of arbitrarily high order. This allows us to view the
representation theory of semi simple Lie algebras from the alternative
perspective proposed in [O1].
\end{abstract}

\section{Jet-filtrations on representations}

Let $\rho :\mathfrak{g\rightarrow gl(}V)$ be a representation of a (finite
dimensional) Lie algebra $\mathfrak{g.}$ Suppose we have an ascending
filtration of subspaces

\begin{equation}
\{0\}\varsubsetneqq V_{1}\varsubsetneqq V_{2}\varsubsetneqq
....\varsubsetneqq V_{k}\varsubsetneqq V
\end{equation}%
satisfying

\begin{equation}
V_{i}\varsubsetneqq \rho (\mathfrak{g})(V_{i})\subset V_{i+1}\text{ \ \ \ \
\ }1\leq i\leq k
\end{equation}%
where $\rho (\mathfrak{g})(V_{i})\overset{def}{=}\{\rho (x)(v)\mid x\in 
\mathfrak{g,}$ $v\in V_{i}\}$ and $V_{k+1}=V$.

\begin{definition}
A filtration (1) satisfying (2) is a jet-filtration on the representation $%
\rho .$ The integer $k$ is the length of the jet filtration (1).
\end{definition}

We can refine (1): For instance, if there exists some $V_{i}\varsubsetneqq
W_{i}\varsubsetneqq V_{i+1}$ satisfying $V_{i}\varsubsetneqq \rho
(V_{i})\subset W_{i}\varsubsetneqq \rho (\mathfrak{g})(W_{i})\subset
V_{i+1}, $ then we can insert $W_{i}$ into (1) and obtain a finer
jet-filtration. In particular, we define a maximally refined jet-filtration
in the obvious way. In the opposite direction, we can omit some term(s) from
(1) and the remaining terms define a coarser jet-filtration.

\begin{definition}
The jet-order of the representation $\rho :\mathfrak{g\rightarrow gl(}V)$ is
the maximum of the lengths of all (maximally refined) jet-filtrations (1).
\end{definition}

Three natural questions:

\bigskip

\textbf{Q1. }What does (1) have to do with jets?

\bigskip

\textbf{Q2. }How do we construct jet-filtrations?

\bigskip

\textbf{Q3. }Why are jet-filtrations relevant?

\bigskip

Let us start with \textbf{Q2. }

\section{Irreducible representations}

Suppose that $\rho :\mathfrak{g\rightarrow gl(}V)$ is irreducible. We choose
some $0\neq v\in V$ and define $V_{1}\overset{def}{=}Span\{v\}.$ If $%
V_{1}=V, $ we stop. If not, then $\rho (\mathfrak{g)}(V_{1})\nsubseteqq
V_{1} $ by irreducibility. In this case, we define $V_{2}\overset{def}{=}%
Span\{\rho (\mathfrak{g)}(V_{1})\cup V_{1}\}$ and get $V_{1}\varsubsetneqq
V_{2}$ and $V_{1}\varsubsetneqq \rho (\mathfrak{g)}(V_{1})\subset V_{2}.$ If 
$V_{2}=V$ we stop. If not, then $\rho (\mathfrak{g)}(V_{2})\nsubseteqq V_{2}$
by irreducibility. In this case, we define $V_{3}\overset{def}{=}Span\{\rho (%
\mathfrak{g)}(V_{2})\cup V_{2}\},$ and get $V_{2}\varsubsetneqq V_{3}$ and $%
V_{2}\varsubsetneqq \rho (\mathfrak{g)}(V_{2})\subset V_{3}.$ Continuing
this process, we finally get a jet-filtration (1) such that $%
V_{k-1}\varsubsetneqq \rho (\mathfrak{g)}(V_{k})=V.$ Note that this
jet-filtration is maximally refined by construction and its length depends
on the vector $v$ we start with.

\begin{definition}
A vector which maximizes the lengths of all such jet-filtrations is called a
maximal vector of the irreducible representation $\rho .$
\end{definition}

The above construction gives an algorithm for constructing jet-filtrations
using an arbitrary irreducible representation. Since irreducible
representations of semi simple Lie algebras are well known, it is
instructive at this stage to look at a concrete example. So let $\mathfrak{g}%
=\mathfrak{sl}(2,\mathbb{R})$ with the well known basis%
\begin{equation}
e=\left[ 
\begin{array}{cc}
0 & 1 \\ 
0 & 0%
\end{array}%
\right] ,\text{ \ \ }f=\left[ 
\begin{array}{cc}
0 & 0 \\ 
1 & 0%
\end{array}%
\right] ,\text{ \ \ }h=\left[ 
\begin{array}{cc}
1 & 0 \\ 
0 & -1%
\end{array}%
\right]
\end{equation}%
and $V=\mathbb{R}_{k}\mathbb{[}x,y]=$ polynomials in the variables $x,y$ of
total degree $\leq k.$ Now $\dim V=k+1$ and has the basis $v_{1}=y^{k},$ $%
v_{2}=y^{k-1}x,...,v_{k}=yx^{k-1},v_{k+1}=x^{k}.$ We define the linear map $%
\rho :\mathfrak{sl}(2,\mathbb{R})\rightarrow End(\mathbb{R}_{k}\mathbb{[}%
x,y])$ by giving its values on this basis as

\begin{equation}
\rho (e)\overset{def}{=}x\frac{\partial }{\partial y},\text{ \ \ }\rho (f)%
\overset{def}{=}y\frac{\partial }{\partial x},\text{ \ \ }\rho (h)\overset{%
def}{=}x\frac{\partial }{\partial x}-y\frac{\partial }{\partial y}
\end{equation}%
and check by a straightforward computation that (4) gives a representation $%
\rho :\mathfrak{sl}(2,\mathbb{R})\rightarrow \mathfrak{gl}(\mathbb{R}_{k}%
\mathbb{[}x,y])$. Now $\rho (e)v_{1}=v_{2},$ $\rho (e)v_{2}=v_{3},...,\rho
(e)v_{k}=v_{k+1}$ and $\rho (e)v_{k+1}=0.$ We define $V_{i}=Span%
\{v_{1},v_{2},...,v_{i})$ and get the filtration

\begin{equation}
\{0\}\varsubsetneqq V_{1}\varsubsetneqq V_{2}\varsubsetneqq
....\varsubsetneqq V_{k}\varsubsetneqq V
\end{equation}%
satisfying

\begin{equation}
V_{i}\varsubsetneqq \rho (e)(V_{i})\subset V_{i+1}
\end{equation}

The key fact now is that both $\rho (f)$ and $\rho (h)$ preserve the
filtration (5) because $\rho (f)v_{i}=v_{i-1}$ and $\rho (h)v_{i}=\lambda
_{i}v_{i}$ for some constant $\lambda _{i}.$ Therefore we can replace $e$ in
(6) by the whole $\mathfrak{sl}(2,\mathbb{R})$ and (5) becomes a
jet-filtration on the irreducible representation $\rho .$ Thus for any
positive integer $k,$ the above well known irreducible representation of $%
\mathfrak{sl}(2,\mathbb{R})$ of degree $k+1$ has the jet-filtration (5) with
jet-order $k$ and has $v_{1}$ as a maximal vector$.$

It is well known that irreducible representations of semi simple Lie
algebras are "made up" of the representations of $\mathfrak{sl}(2,\mathbb{R}%
).$(see, for instance, [H]). Without going into the technical details, we
will state here the following proposition whose proof is now almost trivial
for anyone familiar with the representation theory of semi simple Lie
algebras.

\begin{proposition}
Let $\mathfrak{g}$ be any semi simple Lie algebra and $k$ be any positive
integer. Then there exists an irreducible representation of $\mathfrak{g}$
whose jet-order is greater than $k.$
\end{proposition}

We now come to \textbf{Q1. }

\section{Klein geometries}

Let $(\mathfrak{h,h}_{0})$ be an infinitesimal and effective Klein geometry,
i.e., $\mathfrak{h}$ is a finite dimensional Lie algebra, $\mathfrak{h}%
_{0}\subset \mathfrak{h}$ is a subalgebra which does not contain any ideals
of $\mathfrak{h}$ other than $\{0\}$ (see [O1], [O3] for details). Any such
Klein geometry defines the filtration (called the Weissfeiler filtration in
some works)

\begin{equation}
\mathfrak{h\supsetneqq h}_{0}\supsetneqq \mathfrak{h}_{1}\supsetneqq
....\supsetneqq \mathfrak{h}_{m}\supsetneqq \mathfrak{h}_{m+1}=\{0\}
\end{equation}%
where $\mathfrak{h}_{i+1}\overset{def}{=}\{x\in \mathfrak{h}_{i}\mid \lbrack
x,\mathfrak{h]\subset h}_{i}\},$ $0\leq i\leq k.$ We called the integer $m+1$
the infinitesimal order of $(\mathfrak{h,h}_{0})$ and denoted it by $ord(%
\mathfrak{h,h}_{0}).$ Now let $(G,H)$ be an effective Klein geometry
inducing $(\mathfrak{h,h}_{0}),$ i.e., $G$ is a Lie group with Lie algebra $%
\mathfrak{g,}$ $H\subset G$ a closed subgroup with Lie algebra $\mathfrak{h}%
. $ Now $G$ acts on $G/H=M$ as a transitive transformation group with the
stabilizer $H.$ If $g(x)=y,$ then the transformation $g\in G$ is locally
determined near $x$ (sometimes globally on $M!)$ by its $m$'th order jet at $%
x$ where $m=ord(\mathfrak{h,h}_{0}).$ Therefore, in order to answer \textbf{%
Q1, }we must relate Klein geometries to representations.

Now given some representation $\rho :\mathfrak{g\rightarrow gl(}V),$ we
define the operation $[,]$ on $\mathfrak{g\times }V$ by

\begin{equation}
\lbrack (g_{1},v_{1}),(g_{2},v_{2})]\overset{def}{=}([g_{1},g_{2}],\rho
(g_{1})(v_{2})-\rho (g_{2})(v_{1}))
\end{equation}

Clearly $[,]$ is bilinear and skew symmetric and a straightforward
verification shows that it also satisfies the Jacobi identity. Therefore $(%
\mathfrak{g\times }V,[,])$ is a Lie algebra which we will denote by $%
\mathfrak{g\times }_{\rho }V$. This construction is a special case of a more
general one described in [J] (see pg. 17). We identify $\mathfrak{g}$ with
its image $(\mathfrak{g},0)\subset \mathfrak{g\times }_{\rho }V$ and
identify $V$ with the abelian ideal $(0,V)\subset \mathfrak{g\times }_{\rho
}V.$ From (8) we deduce the important formula

\begin{equation}
\lbrack \mathfrak{g},V]=[(\mathfrak{g},0),(0,V)]=(0,\rho (\mathfrak{g}%
)V)=\rho (\mathfrak{g})V
\end{equation}

Now suppose that $\rho :\mathfrak{g\rightarrow gl(}V)$ is irreducible with
the jet-filtration (5). To be consistent with the notation in (7), we
rewrite the ascending filtration (5) as a descending filtration (writing $%
V_{k-i}$ for $V_{i})$

\begin{equation}
V\supsetneqq V_{0}\supsetneqq V_{1}\supsetneqq ....\supsetneqq
V_{k}\supsetneqq \{0\}
\end{equation}

We define

\begin{equation}
\mathfrak{h}\overset{def}{=}\mathfrak{g\times }_{\rho }V
\end{equation}

\begin{equation}
\mathfrak{h}_{0}\overset{def}{=}V_{0}\varsubsetneqq V\subset \mathfrak{h}
\end{equation}%
and consider the Klein geometry $(\mathfrak{h,h}_{0})$ which is effective by
(9) and irreducibility of $\rho .$ By the definition of (7), we have $[%
\mathfrak{h}_{i+1},\mathfrak{h}]\subset \mathfrak{h}_{i}$ and $\mathfrak{h}%
_{i+1}$ is the largest subalgebra of $\mathfrak{h}_{i}$ satisfying $[%
\mathfrak{h}_{i+1},\mathfrak{h}]\subset \mathfrak{h}_{i}.$ For $(\mathfrak{%
h,h}_{0})$ defined by (11), (12), we have $[V_{i+1},\mathfrak{h]\subset }%
V_{i}$ by (9) and therefore $V_{i}\subset \mathfrak{h}_{i}.$ It follows that 
$m\geq k.$ If (10) is maximally refined, we conclude $V_{i}=\mathfrak{h}_{i}$
and $m=k.$ Hence we obtain

\begin{proposition}
Let $\rho :\mathfrak{g\rightarrow gl(}V)$ be an irreducible representation
with some maximally refined jet-filtration (10) of length $k.$ Then the
Klein geometry $(\mathfrak{h,h}_{0})$ defined by (11), (12) is effective
with $ord(\mathfrak{h,h}_{0})=k+1.$
\end{proposition}

Proposition 5 gives a partial answer to the first fundamental problem of jet
theory (FP1) posed in [O1]. It is easy to show that $[\mathfrak{h}_{i},%
\mathfrak{h}_{j}]\subset \mathfrak{h}_{i+j}$ in (7) and therefore $\mathfrak{%
h}_{i}$ is abelian if $2i\geq m+1=ord(\mathfrak{h,h}_{0}),$ i.e., the " half
tail" of the Weissfeiler filtration always consists of abelian ideals of $%
\mathfrak{h}_{0}.$ In the examples we construct by (11), (12), note that $%
\mathfrak{h}_{0}$ is always abelian! It is easy to modify our construction
by defining $\mathfrak{h}_{0}\overset{def}{=}$ $\mathfrak{r\times }_{\rho
}V_{0}$ (rather than by (12)) for some "carefully chosen" subalgebra $%
\mathfrak{r\subset g}$ which will make $\mathfrak{h}_{0}$ nonabelian (but
solvable!) and still give the conclusion of Proposition 5. However, we do
not know the degree of freedom we have for $\mathfrak{h}_{0}.$

To conclude this section, we we would like to express our belief that the
modern representation theory of Lie algebras can be reconstructed, possibly
with some new insights and perspectives, as a special case of the theory of
Klein geometries.

We now come to \textbf{Q3 }(which is partially answered by the above
paragraph)\textbf{.}

\section{Stiffening of symmetries and integrability}

Any introductory book on differential geometry starts by defining the
tangent space, vector fields, tensor fields, differential forms...etc. These
are all first order objects, i.e., a transformation acts on these objects by
its first order derivatives. One is naturally lead to ask the following
admittedly vague question:

\textbf{Q: }Do the "higher order derivatives " play any role in the global
structure of smooth manifolds?

If the answer to \textbf{Q }is affirmative, then "higher order derivatives"
are surely quite relevant since understanding the structures of manifolds is
a fundamental problem in differential geometry. We do not know of any
attempt in the literature to answer \textbf{Q} other than the striking
unpublished preprint [Ol].

Now suppose some Lie group $G$ (rather, a pseudogroup of finite order) acts
transitively on a smooth manifold $M.$ This situation expresses the fact
that $M$ is "symmetric" in some particular way. If some proper subgroup $%
G^{\prime }\subset G$ also acts transitively on $M,$ i.e., if the action of
the transformation group $G$ can be stiffened to the transformations of $%
G^{\prime },$ we can speculate that $M$ is "more symmetric" in the same way.
As an extreme, if $G^{\prime \prime }\subset G^{\prime }\subset G$ acts
simply transitively, then $M$ is a Lie group, the most symmetric object from
this standpoint. However, if we allow the transformation group to be "too
large" like $Diff(M),$ then all manifolds become symmetric since $Diff(M)$
acts transitively on $M.$

Now the key fact: For some fixed $M,$ if a Lie group $G,$ whose dimension is
much greater than the dimension of $M,$ acts transitively and effectively on 
$M,$ in which case the dimension of the stabilizer $H$ must be large too
since $\dim G/H=$ $\dim G-\dim H$ $=\dim M,$ then the Klein geometry $(G,H)$
must have high order as defined above, i.e., we need higher order
derivatives of the transformations of $G$ to determine them. The reason is
that $H$ injects into the $k$'th order jet group $G_{k}(n)$ where $n=\dim M.$
For instance, if $k=1,$ then we must have $\dim H\leq \dim G_{1}(n)=n^{2}.$
If $\dim H\succneqq \dim G_{1.}(n),$ then we need at least second order
derivatives to inject $H$ into $G_{2}(n),$ if $\dim H\succneqq \dim
G_{2.}(n),$ then we need at least third order derivatives...etc. In short,
"large groups act transitively and effectively on small spaces with large
order". Now suppose $(G,H)$ has order $k+1$ and the action of $G$ can not be
stiffened to a proper subgroup of order $k.$ This indicates an obstruction
coming from $k$'th order derivatives and the problem is, of course, to
formulate this obstruction in precise mathematical terms. Now $(G,H)$
defines a flat pre-homogeneous geometry (PHG) of order $k+1$ ([O1], [O2])
which will "contain" many PHG's of order $k$ and none of these PHG's can
integrate, i.e., have vanishing curvature, for otherwise they would give a
stiffening.

Remarkably, the stiffening problem has a purely algebraic counterpart: If $%
(G,H)$ stiffens $(G^{\prime },H^{\prime }),$ then $(\mathfrak{%
g,h)\preccurlyeq (g}^{\prime },\mathfrak{h}^{\prime }),$ i.e., $\mathfrak{g+h%
}^{\prime }=\mathfrak{g}^{\prime },$ $\mathfrak{h=g\cap h}^{\prime }$ (see
[O1], pg. 131 for details). Note that $\dim \mathfrak{g-}\dim \mathfrak{h=}%
\dim \mathfrak{g}^{\prime }-\dim \mathfrak{h}^{\prime }.$ Now $\preccurlyeq $
is a partial order on Klein geometries and the problem is to find the
minimal (also maximal!) elements.

The above speculative scenario has been one of the motivations for us to
become interested in Klein geometries and geometric structures of high order
and we hope it may be inspiring also for some others.

\bigskip

\textbf{References}

\bigskip

[H] B.C.Hall: Lie Groups, Lie Algebra, and Representations: An Elementary
Introduction, Graduate Texts in Mathematics, Book 222, Springer, 2015

[J] N. Jacobson: Lie Algebras, Dover Publications Inc., 1962

[Ol] P.J.Olver: Differential hyperforms, I, preprint, University of
Minnesota, 1982

[O1] E.H.Orta\c{c}gil: An Alternative Approach to Lie Groups and Geometric
structures, OUP, 2018

[O2] \ \ \ \ \ \ " \ \ \ \ \ \ \ \ \ \ : Curvature without connection,
arXiv:2003.06593, 2020

[O3] \ \ \ \ \ "\ \ \ \ \ \ \ \ \ \ \ \ : Nilpotency and higher order
derivatives in differential geometry, arXiv:2105.07254, 2021

\bigskip

Erc\"{u}ment H. Orta\c{c}gil

ortacgile@gmail.com

\end{document}